\documentclass[12pt]{amsart}
\advance\hoffset by -0.5 truecm
\usepackage{amssymb}
\newtheorem{Theorem}{Theorem}[section]
\newtheorem{Lemma}[Theorem]{Lemma}
\newtheorem{Corollary}[Theorem]{Corollary}

\newtheorem{Remark}[Theorem]{Remark}

\def\V{\mbox{Var}}

\def\R\re
\def\V{\bf V}

\def \la{\lambda}

\def \re{{\mathbb R}}

\def \M{{\widetilde M}}

\def \0{\lambda_{0}}
\def \la{\lambda}
\def \ga{\gamma}

\def\wO{\widetilde\Omega}

\begin{document}
\title[Magnetic rigidity of horocycle flows]{Magnetic rigidity of horocycle flows}

\author[G. P. Paternain]{Gabriel P. Paternain}

 \address{ Department of Pure Mathematics and Mathematical Statistics,
University of Cambridge,
Cambridge CB3 0WB, England}
 \email {g.p.paternain@dpmms.cam.ac.uk}

\date{September 2004}


\begin{abstract} Let $M$ be a closed oriented surface endowed with a Riemannian metric $g$ 
and let $\Omega$ be a 2-form. We show that the magnetic flow of the pair
$(g,\Omega)$ has zero asymptotic Maslov index and zero Liouville action if and only $g$ has constant
Gaussian curvature, $\Omega$ is a constant multiple of the area form of $g$ and the magnetic flow
is a horocycle flow.

This characterization of horocycle flows implies that if the magnetic flow of a pair 
$(g,\Omega)$ is $C^1$-conjugate to the horocycle flow of a hyperbolic metric $\bar{g}$
then there exists a constant $a>0$, such that $ag$ and $\bar{g}$ are isometric and 
$a^{-1}\Omega$ is, up to a sign, the area form of $g$.
The characterization also implies that if a magnetic flow is Ma\~n\'e critical and
uniquely ergodic it must be the horocycle flow.

As a by-product we also obtain results on existence of closed magnetic geodesics
for almost all energy levels in the case weakly exact magnetic fields on arbitrary manifolds.

\end{abstract}

\maketitle

\section{Introduction}

Let $\Gamma$ be a cocompact lattice of $PSL(2,\re)$.
The standard horocycle flow $h$ is given by the right action of the one-parameter subgroup
\[\left(\begin{array}{cc}

1&t\\
0&1\\

\end{array}\right)   \]
on $\Gamma\backslash PSL(2,\re)$.
The horocycle flow is known to display very peculiar ergodic properties.
It preserves the Riemannian volume on $\Gamma\backslash PSL(2,\re)$, is
uniquely ergodic \cite{Fu}, and mixing of all degrees \cite{M1}. Morever, it has zero entropy
since 
\begin{equation}
\phi^0_t\circ h_s=h_{se^{-t}}\circ \phi^0_t
\label{con1}
\end{equation} 
for all $s,t\in \re$, where $\phi^0$ is the geodesic flow given by the one-parameter subgroup:
\[\left(\begin{array}{cc}

e^{t/2}&0\\
0&e^{-t/2}\\

\end{array}\right).   \]
In fact, $h$ parametrizes the strong stable manifold of $\phi^0$.
The horocycle flow is a very rigid object, as the works of B. Marcus and M. Ratner show
\cite{M,Ra1,Ra2}. Recent results on ergodic averages and solutions of cohomological equations 
for $h$ can be found in \cite{Bu,FF}.

In the present paper we would like to look at horocycle flows as magnetic flows.
A matrix $X$ in $sl(2,\re)$ gives rise to a flow $\phi$ on $\Gamma\backslash PSL(2,\re)$
by setting
\[\phi_{t}(\Gamma g)=\Gamma g e^{tX}.\]
The geodesic and horocycle flows are just particular cases of these algebraic flows.
Consider the following path of matrices in $sl(2,\re)$:
\[\re\ni\la\mapsto X_{\la}:=\left(\begin{array}{cc}

1/2&0\\
0&-1/2\\ 
\end{array}\right)+\la\,\left(\begin{array}{cc}

0&1/2\\
-1/2&0\\

\end{array}\right).  \]
The flows $\phi^{\la}$ on $\Gamma\backslash PSL(2,\re)$ associated with the matrices $X_{\la}$
have an interesting interpretation. Since $PSL(2,\re)$ acts by isometries on the hyperbolic
plane ${\mathbb H}^{2}$, $M:=\Gamma\backslash {\mathbb H}^{2}$ is a compact
hyperbolic surface (provided $\Gamma$ acts without fixed points) and the unit sphere bundle $SM$ of $M$ can be identified
with $\Gamma\backslash PSL(2,\re)$. 
A calculation shows that $\phi^{\la}$ is the
Hamiltonian flow of the Hamiltonian
$H(x,v)=\frac{1}{2}|v|^{2}_{x}$ with respect to the symplectic
form on $TM$ given by
\[-d\alpha+\la\,\pi^{*}\Omega_{a},\]
where $\Omega_{a}$ is the area form of $M$, $\pi:TM\to M$ is
the canonical projection and $\alpha$ is the contact 1-form that generates the geodesic 
flow of $M$.
For $\la=0$, $\phi^{0}$ is the geodesic flow and for $\la=1$, $\phi^{1}$ is the flow
induced by the one-parameter subgroup with matrix on $sl(2,\re)$ given by
\[X_1=\left(\begin{array}{cc}

1/2&1/2\\
-1/2&-1/2\\

\end{array}\right). \]
Since there exists an element $c\in PSL(2,\re)$ such that
\[ c^{-1}X_1 c=\left(\begin{array}{cc}
0&1\\
0& 0\\
\end{array}\right) \]
the map $f:\Gamma\backslash PSL(2,\re)\to \Gamma\backslash PSL(2,\re)$ given by $f(\Gamma g)=\Gamma gc$ conjugates
$\phi^1$ and $h$, i.e. $f\circ \phi^{1}_{t}=h_{t}\circ f$. In fact, any matrix
in $sl(2,\re)$ with determinant zero will give rise to a flow which is conjugate
to $h_t$ or $h_{-t}$. (So, up to orientation, there is just one algebraic horocycle
flow.) Passing by, we note that $\det X_{\la}=-\frac{1}{4}(1-\la^2)$, so for $|\la|<1$,
the flow $\phi^{\la}$ is conjugate to the geodesic flow $\phi^0$, up to a constant time
scaling by $\sqrt{1-\la^2}$.
Hence the magnetic flows $\phi^\la$ are just geodesic
flows, but with entropy $\sqrt{1-\la^2}$. (This latter observation is due to 
V.I. Arnold \cite{A1/2}.) 

In general, if $(M,g)$ is a closed Riemannian manifold and $\Omega$ is a closed 2-form,
the Hamiltonian flow $\phi$ of
$H(x,v)=\frac{1}{2}|v|^{2}_{x}$ with respect to the symplectic
form on $TM$ given by
\[\omega:=-d\alpha+\pi^{*}\Omega,\]
is called the {\it magnetic flow} of the
pair $(g,\Omega)$ because it models the motion of a particle under the influence of
the magnetic field $\Omega$. The projection of the orbits of $\phi$ to $M$ are
called {\it magnetic geodesics}. The discussion above shows that $h$ appears as the
magnetic flow of a hyperbolic surface with $\Omega$ equal to the area form of the surface.

Since the horocycle flow has no closed orbits, this already gives an example
of a Hamiltonian system with an energy level ($SM$) without closed orbits. This example
has been much used, most notably by V. Ginzburg \cite{Gi1,Gi2,Gi3} to give smooth counterexamples
to the Hamiltonian Seifert conjecture in $\re^{2n}$, $n\geq 3$ 
(a $C^2$-counterexample is now available in $\re^4$ \cite{GG}).

How frequently does the horocycle flow appear as a magnetic flow? To answer this question
we first prove a characterization of horocycle flows within the set of magnetic flows.

Magnetic flows on surfaces leave invariant the volume form $\alpha\wedge d\alpha$.
The associated Borel probability measure is called the {\it Liouville measure} $\mu_{\ell}$
of $SM$. We shall assume from now on that
$M$ has genus $\geq 2$. Then $\pi^{*}:H^{2}(M,\re)\to H^{2}(SM,\re)$
is the zero map and thus if $\Omega$ is any 2-form, $\pi^{*}\Omega$ is exact on $SM$.
It follows that $\omega$ restricted to $SM$ is exact and we let $\Theta$ be any primitive.
We define {\it the action of the Liouville measure} as:
\[{\frak a}(\mu_{\ell}):=\int\Theta(X)\,d\mu_{\ell},\]
where $X$ is the vector field on $SM$ that generates the magnetic flow. 
The action does not depend on the primitive, since the {\it asymptotic cycle} of 
$\mu_{\ell}$ is zero (cf. Section 2 and \cite{CGIP}), i.e. for any closed 1-form $\varphi$ on $SM$ we have
\[\int\varphi(X)\,d\mu_{\ell}=0.\]
It is quite simple to check that when $M$ is a hyperbolic surface and $\Omega$
is the area form, ${\frak a}(\mu_{\ell})=0$. In fact, in this case, there is
a primitive $\Theta$, with $\Theta(X)\equiv 0$. 
It is also easy to check that there are
no conjugate points \cite{CI}[Example A.1]. Equivalently, using the results
in \cite{CGIP}, we can say that {\it the asymptotic Maslov index} ${\frak m}(\mu_{\ell})$
of the Liouville measure is zero (cf. Section 2).

We first show that these two symplectic-ergodic quantities characterize horocycle flows.

\medskip

\noindent {\bf Proposition.} {\it Let $M$ be a closed oriented surface endowed with
a Riemannian metric $g$ 
and let $\Omega$ be a 2-form. The magnetic flow of the pair
$(g,\Omega)$ has ${\frak a}(\mu_{\ell})={\frak m}(\mu_{\ell})=0$ if and only 
$g$ has constant Gaussian curvature, $\Omega$ is a constant multiple of the area form of $g$ 
and the magnetic flow is a horocycle flow.

}

\medskip

The Proposition has the following consequence:

\medskip

\noindent {\bf Theorem A.} {\it Let $M$ be a closed oriented surface endowed with
a Riemannian metric $g$ 
and let $\Omega$ be a 2-form. If the magnetic flow of the pair 
$(g,\Omega)$ is $C^1$-conjugate to the horocycle flow of a hyperbolic metric $\bar{g}$,
there exists a constant $a>0$, such that $ag$ and $\bar{g}$ are isometric and 
$a^{-1}\Omega$ is, up to a sign, the area form of $g$.
}

\medskip

We observe that it is not possible to conclude that $g$ has curvature $-1$ (i.e. $a=1$).
This is simply because the magnetic flow of a pair $(g,\Omega)$ with $g$ of constant
negative curvature $-k$, $\Omega=\la\,\Omega_a$ ($\la>0$) and $\la^2=k$ is smoothly
conjugate to the horocyle flow of a hyperbolic surface. To see this observe that
an easy scaling argument shows that the magnetic flow of a such pair is, up to a constant
time change, conjugate to the horocycle flow of a hyperbolic surface.
But by (\ref{con1}), $h_t$ is conjugate to $h_{\tau t}$ for any positive real number
$\tau$.
This shows that the area $A$ of a surface is {\it not} preserved under $C^1$-conjugacies
of magnetic flows. C. Croke and B. Kleiner \cite{CK} have shown that the volume
of a Riemannian manifold is preserved under $C^1$-conjugacies of geodesic flows.
In the case of transitive magnetic flows, we show that $C^1$-conjugacies
preserve ${\mathfrak a}(\mu_{\ell})/A$ (cf. Lemma \ref{l1}).

We now describe a second and more involved application of the Proposition.
Let $\wO$ be the lift of $\Omega$ to  the universal 
cover~$\widetilde M\cong\re^2$ of~$M$. 
Since $\widetilde\Omega$ is an exact form, there
exists a smooth 1-form $\theta$ such that $\wO=d\theta$.
Let us consider the Lagrangian on $\M$ given by
\[L(x,v)=\frac{1}{2}|v|^{2}_{x}-\theta_{x}(v).\]
It is well known that the extremals of $L$, i.e., the solutions of
the Euler-Lagrange equations of $L$,
 $$
 \frac{d}{dt}\frac{\partial L}{\partial v}(x,v) = \frac{\partial L}{\partial x}(x,v)
$$
coincide with the lift to $\M$ of the
magnetic geodesics. 
The action of the Lagrangian $L$ on
an absolutely continuous curve $\ga:[a,b]\rightarrow \M$ is defined by
\[A_{L}(\ga)=\int_{a}^{b}L(\ga(t),\dot{\ga}(t))\,dt.\]
The {\it Ma\~n\'e critical value} of the pair $(g,\Omega)$ is
\[c=c(g,\Omega):=\inf\{k\in\re:\;A_{L+k}(\ga)\geq 0\;\mbox{\rm  
for any absolutely continuous closed curve $\ga$}\]
\[\mbox{\rm defined on any closed interval $[0,T]$ }\}.\]
Like any Lagrangian flow, the magnetic flow for $T\M$ can be viewed as the Hamiltonian flow 
defined by the canonical symplectic form on $T^*\M$ and a 
suitable Hamiltonian function $H:T^*\M \to \re$; in this case
\[ H(x,p) = \frac12 |p + \theta_x|^2. \]
The Legendre transform ${\mathcal L}:T\M\to T^{*}\M$ defined
by
$${\mathcal L}(x,v)=\frac{\partial L}{\partial v}(x,v)$$
carries orbits of the Lagrangian flow for $L$ to orbits of the Hamiltonian flow defined
by $H$ and the canonical symplectic form.
The critical value can also be defined in Hamiltonian terms \cite{BP} as:
\begin{align*}
c(g,\Omega) &= \inf_{u\in C^{\infty}(\M,\re)}\;\sup_{x\in
  \M}\;H(x,d_{x}u)\\
&= \inf_{u\in C^{\infty}(\M,\re)}\;\sup_{x\in \M}\;
   \frac{1}{2}|d_{x}u+\theta_{x}|^{2}.
\end{align*}
As $u$ ranges over $C^{\infty}(\M,\re)$ the form $\theta-du$ ranges over all
primitives of $\wO$, because any two primitives differ by a closed 1-form
which must be exact since $\M$ is simply connected.
Since on a surface of genus $\geq 2$, there are bounded primitives we always have
$c(g,\Omega)<\infty$.

We will say that a magnetic flow is {\it Ma\~n\'e critical} if $c(g,\Omega)=1/2$.
Magnetic flows which are supercritical, i.e. $1/2>c(g,\Omega)$ always have positive 
topological entropy \cite{BP}[Proposition 5.4]. Hence such flows exhibit a horseshoe
and exponential growth rate of hyperbolic closed magnetic geodesics. In fact, we will show
that, in any dimension, if $1/2>c(g,\Omega)$, then a nontrivial homotopy class contains
a closed magnetic geodesic provided that the centralizer of some element
in the class is an amenable subgroup (cf. Theorem \ref{cloorb}). 
For subcritical magnetic flows, i.e. $1/2<c(g,\Omega)$ one hopes to prove that there
are always closed contractible magnetic geodesics, although nothing of this kind has been
proved in general. When $\Omega$ itself is exact, the main result in
\cite{CMP} says that there always exists a closed magnetic geodesic.

What happens for magnetic flows which are Ma\~n\'e critical with $\Omega$ non-exact?
It is easy to check that the hororcycle flow is Ma\~n\'e critical (cf. \cite{Co}[Example  6.2]). Is it the only 
magnetic flow which is Ma\~n\'e critical and uniquely ergodic?
Aubry-Mather theory combined with the Proposition gives the answer:

\medskip

\noindent {\bf Theorem B.} {\it Let $M$ be a closed oriented surface endowed with
a Riemannian metric $g$ 
and let $\Omega$ be a 2-form. Suppose the magnetic flow of the pair 
$(g,\Omega)$ is Ma\~n\'e critical and uniquely ergodic.
Then $g$ has constant Gaussian curvature, 
$\Omega$ is a constant multiple of the area form of $g$ 
and the magnetic flow is a horocycle flow.
}

\medskip

Our results on the existence of closed magnetic geodesics in non-trivial free
homotopy classes for supercritical magnetic flows combined with recent
results of G. Contreras \cite{Co1} and O. Osuna \cite{O} imply the following statement
on almost existence of closed magnetic geodesics for weakly exact magnetic flows.
Recall that $\Omega$ is said to be {\it weakly exact} if its lift to the universal
covering of $M$ is exact.

\medskip

\noindent {\bf Theorem C.} {\it Let $M$ be an arbitrary closed manifold endowed with
a Riemannian metric $g$ and let $\Omega$ be a weakly exact 2-form. 
We have:
\begin{enumerate}
\item if $\pi_1(M)$ is amenable and $\Omega$ is not exact, then almost every energy
level contains a contractible closed magnetic geodesic;
\item if $\pi_1(M)$ contains a non-trivial element with an amenable centralizer, then
almost every energy level contains a closed magnetic geodesic.
\end{enumerate}
}

\medskip

Recall that a discrete group $\Gamma$ is said to be {\it amenable} if the space of bounded functions
$\Gamma\to\re$ has a left (or right) invariant mean \cite{Pi}. Examples are finite groups, abelian groups
and finite extensions of solvable groups. If a group contains a free subgroup on two generators, then is non-amenable.
I do not know of any example of a finitely presented group, for which the centralizer
of every element is non-amenable. We note that Contreras in \cite{Co1} proved item 2 in Theorem C
when $\Omega$ is exact without any assumption on $\pi_1(M)$ and V. Ginzburg and E. Kerman proved
item 2 when $M$ is a torus \cite{GK}.

There are several results establishing the existence of closed contractible magnetic
geodesics for almost every {\it low} energy level, see \cite{S,Mac} and references therein. 
The methods of Symplectic Topology have proven to be effective in this respect. For high energy levels,
at least when $g$ is generic, one can also obtain existence of closed magnetic geodesics
in free homotopy classes by observing that magnetic flows approach geodesics flows
as energy increases. Very little is known about how to bridge
the gap. The exceptions seem to be Theorem C and the main result in \cite{CMP} alluded above,
which are both based on Aubry-Mather theory.

\medskip

{\it Acknowledgements:} I would like to thank Leonardo Macarini for several useful
comments and discussions about Theorem C. He suggested 
the possibilty of using Contreras' recent results to prove statements like
item 2 in Theorem C. I also thank Viktor Ginzburg for several discussions regarding
the question of existence of closed magnetic geodesics.

\section{Geometric preliminaries}

Let $M$ be a closed oriented surface, $SM$ the unit sphere bundle
and $\pi:SM\to M$ the canonical projection. The latter is in fact
a principal $S^{1}$-fibration and we let $V$ be the infinitesimal
generator of the action of $S^1$.

Given a unit vector $v\in T_{x}M$, we will denote by $iv$ the
unique unit vector orthogonal to $v$ such that $\{v,iv\}$ is an
oriented basis of $T_{x}M$. There are two basic 1-forms $\alpha$
and $\beta$ on $SM$ which are defined by the formulas:
\[\alpha_{(x,v)}(\xi):=\langle d_{(x,v)}\pi(\xi),v\rangle;\]
\[\beta_{(x,v)}(\xi):=\langle d_{(x,v)}\pi(\xi),iv\rangle.\]
The form $\alpha$ is precisely the contact form that we mentioned
in the introduction. The vector field $X_0$ uniquely determined by
the equations $\alpha(X_{0})\equiv 1$, $i_{X_{0}}d\alpha\equiv 0$
generates the geodesic flow $\phi^0$ of $M$. 

A basic theorem in 2-dimensional Riemannian geometry asserts that
there exists a unique 1-form $\psi$ on $SM$ (the connection form)
such that $\psi(V)\equiv 1$ and

\begin{align*}
&d\alpha =\psi\wedge \beta\\ & d\beta=-\psi\wedge
\alpha\\ & d\psi=-(K\circ\pi)\,\alpha\wedge\beta
\end{align*}
where $K$ is the Gaussian curvature of $M$. In fact, the form
$\psi$ is given by
\[\psi_{(x,v)}(\xi)=\left\langle \frac{DZ}{dt}(0),iv\right\rangle,\]
where $Z:(-\varepsilon,\epsilon)\to SM$ is any curve with
$Z(0)=(x,v)$ and $\dot{Z}(0)=\xi$ and $\frac{DZ}{dt}$ is the
covariant derivative of $Z$ along the curve $\pi\circ Z$.

It is easy to check that $\alpha\wedge\beta=\pi^{*}\Omega_{a}$,
hence
\begin{equation}
d\psi=-\pi^{*}(K\,\Omega_{a}). \label{psi}
\end{equation}

In the case of a hyperbolic surface, the vertical vector field $V$
corresponds to the matrix in $sl(2,\re)$ given by
\[\left(\begin{array}{cc}

0&1/2\\
-1/2&0\\

\end{array}\right).  \]

\subsection{Asymptotic cycles} 
Given a $\phi$-invariant Borel probability measure $\mu$, {\it the asymptotic cycle}
of $\mu$ is the real 1-homology class ${\mathcal S}(\mu)$ defined by the equality:
\[\langle [\varphi],{\mathcal S}(\mu)\rangle =\int_{SM}\varphi(X)\,d\mu\]
for any closed 1-form $\varphi$.
Let us check that ${\mathcal S}(\mu_{\ell})=0$. We will prove something slightly more general
that we will need later. Let $\Lambda_{0}^*$ be the space of continuous forms $\lambda$
whose exterior derivative, weakly defined by Stokes' theorem: 
$\int_{\sigma} d\la=\int_{\partial \sigma}\la$ for every smooth chain $\sigma$,
are also continuous differential forms. The space $\Lambda_{0}^*$ is closed
under exterior differentiation, wedge products and pull back of $C^1$ maps.

Let $\varphi$ be a continuous 1-form
which is closed in the sense that its integral over the boundary of any 2-chain
is zero, i.e. $\varphi\in \Lambda_{0}^1$ and $d\varphi=0$. We claim that
\[ \int_{SM}\varphi(X)\,d\mu_{\ell}=0.\]
Note that we always have $\alpha\wedge \pi^{*}\Omega=0$ as it easily follows
from evaluating the 3-form on any basis that contains $V$. Thus
$\alpha\wedge(-d\alpha)=\alpha\wedge \omega$ and it suffices to show that
\[ \int_{SM}\varphi(X)\,\alpha\wedge \omega=0.\]
Observe that 
$$\varphi(X)\,\alpha\wedge\omega=\varphi\wedge i_{X}(\alpha\wedge\omega)=\varphi\wedge\omega,$$
since $i_{X}\omega=0$. But since $\omega$ is exact, if we let $\Theta$ be a primitive, we have
$d(\varphi\wedge\Theta)=\varphi\wedge\omega$ and the claim follows from our definition
of exterior differentiation and the fact that $M$ is a closed surface.

\subsection{Asymptotic Maslov index}

Let $\Lambda(SM)$ be the set of Lagrangian subspaces of $T|_{SM}TM$.
Given any subspace $E\in \Lambda(SM)$ and $T>0$ we can consider the curve
of Lagrangian subspaces $[0,T]\ni t\mapsto d\phi_{t}(E)$. Let
$n(x,v,E,T)$ be the intersection number of this curve with the Maslov cycle
of $\Lambda(SM)$. It was shown in \cite{CGIP} that if $\mu$ is $\phi$-invariant,
the limit
\[{\mathfrak m}(x,v):=\lim_{T\to\infty}\frac{1}{T}n(x,v,E,T)\]
exists for $\mu$-almost every $(x,v)$, is independent of $E$, and $(x,v)\mapsto {\mathfrak m}(x,v)$
is integrable. {\it The asymptotic Maslov index} of $\mu$ is:
\[{\mathfrak m}(\mu):=\int_{SM}{\mathfrak m}(x,v)\,d\mu(x,v).\]

\subsection{Green subbundles}

Given $(x,v)\in TM$ we define {\it the vertical subspace
at} $(x,v)$ as ${\mathcal V}(x,v):=\mbox{\rm ker}\,d_{(x,v)}\pi$, where $\pi:TM\to M$
is the canonical projection. Note that ${\mathcal V}(x,v)\cap T_{(x,v)}SM$ is spanned
by the value of the vector field $V$ at $(x,v)$.
We say that the orbit of $(x,v)\in SM$ {\it does not have conjugate
points} if for all $t\neq 0$,
$$d_{(x,v)}\phi_{t}({\mathcal V}(x,v))\cap {\mathcal V}(\phi_{t}(x,v))=\{0\}.$$

Since magnetic flows are optical, the main result in \cite{CGIP} says that
$SM$ has no conjugate points (i.e. for all $(x,v)\in SM$, the orbit of
$(x,v)$ does not have conjugate points) if and only if ${\mathfrak m}(\mu_{\ell})=0$.

If $SM$ has no conjugate points, one can construct the so
called {\it Green subbundles} \cite[Proposition A]{CI} given by:
\[E(x,v):=\lim_{t\to +\infty}d\phi_{-t}({\mathcal V}(\phi_{t}(x,v))),\]
\[F(x,v):=\lim_{t\to +\infty}d\phi_{t}({\mathcal V}(\phi_{-t}(x,v))).\] These subbundles are Lagrangian, they never intersect the vertical subspace
and, crucial for us, they are contained in $T(SM)$. Moreover, they vary measurably
with $(x,v)$ and they contain the vector field $X$.

\section{Proof of the Proposition} The proof will be based on integrating an 
appropriate Riccati equation along a solution arising from the Green bundles.
This is a well known method, first used by E. Hopf \cite{Hopf} and subsequently
extended to higher dimensions by L.W. Green \cite{Green}. The method is still
paying dividends, cf. \cite{Bi,Gouda}.

Let $\Omega$ be an arbitrary smooth
2-form. We write $\Omega=f\,\Omega_{a}$, where $f:M\to\re$ is a
smooth function and $\Omega_a$ is the area form of $g$. Let $A$
be the area of $g$.

 Since $H^{2}(M,\re)=\re$, there exist a constant $c$ and a
smooth 1-form $\varrho$ such that
\[\Omega=cK\,\Omega_{a}+d\varrho \]
and $c=0$ if and only if $\Omega$ is exact. Using (\ref{psi}) we
have
\[\omega:=-d\alpha+\,\pi^{*}\Omega=d(-\alpha-
c\,\psi+\pi^{*}\varrho).\] The vector field $X$ that
generates the magnetic flow $\phi$ is given by $X=X_{0}+f\,V$ since it
satisfies the equation $dH=i_{X}\omega$. 
Since $X_0$ and $V$ preserve the volume form $\alpha\wedge d\alpha$,
then so does $X=X_{0}+f V$ and thus $\phi$ preserves
the normalized Liouville measure $\mu_{\ell}$ of $SM$. 

If we evaluate the primitive $-\alpha-c\,\psi+\pi^{*}\varrho$ of
the symplectic form $\omega$ on $X$ we obtain:
\begin{equation}
(-\alpha-c\,\psi+\pi^{*}\varrho)(X)(x,v)=-1-
c\,f(x)+\varrho_{x}(v).\label{contact}
\end{equation}
Therefore ${\frak a}(\mu_{\ell})$ is given by:
\[{\frak a}(\mu_{\ell})=-1-\frac{c}{A}\int _{M}f(x)\,dx\]
since $\mu_{\ell}$ is invariant under the flip $v\mapsto
-v$. By the definition of $c$ and $f$ and the Gauss-Bonnet theorem
\[\int_{M}f(x)\,dx=c\int_{M}K(x)\,dx=2\pi\,\chi\,c\]
and hence
\begin{equation}
{\frak a}(\mu_{\ell})=-1-\frac{1}{2\pi\chi\,A}\left(\int_{M} f(x)\,dx\right)^{2}.
\label{action}
\end{equation}  

Given $(x,v)\in SM$ and $\xi\in T_{(x,v)}TM$, let
\[J_{\xi}(t)=d_{(x,v)}(\pi\circ\phi_{t})(\xi).\]
We call $J_{\xi}$ a {\it magnetic Jacobi field} with initial
condition $\xi$. It was shown in \cite{PP1/2} that $J_{\xi}$
satisfies the following Jacobi equation:
\begin{equation}
\ddot{J_{\xi}} + R(\dot{\gamma},J_{\xi})\dot{\gamma} -[Y(\dot{J_{\xi}})+
(\nabla_{J_{\xi}}Y)(\dot{\gamma})] = 0,  \label{jacobi}
\end{equation}
where $\ga(t)=\pi\circ\phi_{t}(x,v)$, $R$ is the curvature tensor
of $g$ and $Y$ is determined by the equality $\Omega_{x}(u,v)=\langle Y_x(u),v\rangle$
for all $u,v\in T_{x}M$ and all $x\in M$.

Let us express $J_{\xi}$ as follows:
\[J_{\xi}(t)=x(t)\dot{\ga}(t)+y(t)i\dot{\ga}(t),\]
and suppose in addition that $\xi\in T_{(x,v)}SM$, which implies
\begin{equation}
g_{\ga}(\dot{J}_{\xi},\dot{\ga})=0.\label{tangente}
\end{equation}
A straightforward computation using (\ref{jacobi}) and (\ref{tangente}) shows that
$x$ and $y$ must satisfy the scalar equations:

\begin{align}
&\dot{x}=f(\ga)\,y   \label{magJac1}\\
&\ddot{y}+\left[K(\ga)-\left\langle\nabla
f(\ga),i\dot{\ga}\right\rangle+f^{2}(\ga)\right]y=0.
\label{magJac2}
\end{align}

Note that the no conjugate points condition is equivalent to saying that any nontrivial
magnetic Jacobi field which vanishes at $t=0$ is never zero again.

Let us consider one of the Green subbundles, let us say $E$. 
Since for any $(x,v)\in SM$ the subspace $E$ does not intersect
the vertical subspace ${\mathcal V}(x,v)$, there exists a linear map $S(x,v):T_{x}M\to T_{x}M$
such that $E$ can be identified with the graph of $S$.
Let $u(x,v)$ be the trace of $S(x,v)$. An easy calculation using (\ref{magJac1}) and
(\ref{magJac2}) shows that $u$ along $\phi$ satisfies the Riccati equation:
\begin{equation}
\dot{u}+u^{2}+K(\ga)-\left\langle\nabla
f(\ga),i\dot{\ga}\right\rangle+f^{2}(\ga)=0.
\label{riccati}
\end{equation}
We can now integrate equation (\ref{riccati}) with respect to $t\in [0,1]$ and then
with respect to $\mu_{\ell}$ (using the $\phi$-invariance of $\mu_{\ell}$) to conclude that:
\[\int _{SM}u^{2}\,d\mu_{\ell}+\int_{SM}[K(x)-\left\langle\nabla
f(x),i v\right\rangle+f^{2}(x)]\,\mu_{\ell}=0.\]
Since $\mu_{\ell}$ is invariant under the flip $v\mapsto
-v$ we have:
\[\int_{SM}\left\langle\nabla
f(x),i v\right\rangle\,d\mu_{\ell}=0\]
and thus
\[\int _{SM}u^{2}\,d\mu_{\ell}+\int_{SM}[K(x)+f^{2}(x)]\,\mu_{\ell}=0.\]
The last equality implies
\begin{equation}
\int_{SM}[K(x)+f^{2}(x)]\,\mu_{\ell}=\frac{2\pi\chi}{A}+\frac{1}{A}\int_{M}f^{2}(x)\,dx\leq 0
\label{ine1}
\end{equation}
with equality if and only if $u$ is zero for almost every $(x,v)\in SM$.
But if we now assume that the action ${\frak a}(\mu_{\ell})$ vanishes, equation
(\ref{action}) and the Cauchy-Schwarz inequality tell us that:
\[-2\pi\chi\,A=\left(\int_{M} f(x)\,dx\right)^{2}\leq A\,\int_{M}f^{2}(x)\,dx.\]
Combining the last inequality with (\ref{ine1}) we see that
$f$ must be constant and $u$ is zero for almost every $(x,v)\in SM$.
If we now use this information in the Riccati equation (\ref{riccati})
we conclude that $K$ must be constant and $K+f^{2}=0$. The last equality
ensures that the magnetic flow is a horocycle flow thus concluding the proof
of the Proposition.

\section{Proof of Theorem A}

\begin{Lemma} Let $M_i$, $i=1,2$ be closed oriented surfaces
with magnetic flows $\phi^{i}$ determined by pairs $(g_i,\Omega_i)$, $i=1,2$.
Suppose $\phi^1$ is $C^1$-conjugate to $\phi^2$ and one of them
is transitive. Then 
$$A_2\,{\frak a}(\mu^{1}_{\ell})=A_1\,{\frak a}(\mu^{2}_{\ell}),$$
where $A_i$ is the area of $g_i$.
\label{l1}
\end{Lemma}

\begin{proof} Let $f:SM_{1}\to SM_{2}$ be the $C^1$-conjugacy and $\omega_i$
the corresponding symplectic forms restricted to $SM_i$.
Recall that $\alpha\wedge(-d\alpha)=\alpha\wedge \omega$.
Since magnetic flows preserve $\alpha\wedge d\alpha$, the
volume form $f^{*}(\alpha_{2}\wedge\omega_2)$ is invariant under $\phi^1$. 
Since we are assuming that the magnetic flows
are transitive there exists a (nonzero) constant $\kappa$ such that
\begin{equation}
 f^{*}(\alpha_{2}\wedge\omega_2)=\kappa\,\alpha_{1}\wedge\omega_1.
\label{vol}
\end{equation}
Note that $df$ maps $X_1$ to $X_2$, $\alpha_{i}(X_i)=1$ and $i_{X_{i}}\omega_{i}=0$, hence 
contracting with $X_1$ in the last equation gives: \[f^*\omega_2=\kappa\,\omega_1.\]
Let $\Theta_i$ be a primitive of $\omega_i$. Then $\varphi:=f^{*}\Theta_2-\kappa\,\Theta_1$
is a continuous 1-form, which is closed in the sense that its integral
over the boundary of every 2-chain is zero.
By (\ref{vol}), $f_{*}\mu^{1}_{\ell}=\mu^{2}_{\ell}$, thus
\[\int_{SM_{1}}\varphi(X_{1})\,d\mu^{1}_{\ell}=
\int_{SM_{2}}\Theta_{2}(X_{2})\,d\mu^{2}_{\ell}-\kappa\,
\int_{SM_{1}}\Theta_{1}(X_{1})\,d\mu^{1}_{\ell}=
{\frak a}(\mu^{2}_{\ell})-\kappa\,{\frak a}(\mu^{1}_{\ell}).\]
But, since the asymptotic cycle of $\mu_{\ell}$ is zero (cf. Subsection 2.1), the left hand side
vanishes. Equality (\ref{vol}) implies that $\kappa=A_{2}/A_{1}$ and the lemma
follows.

\end{proof}

\begin{Lemma} Let $M_i$, $i=1,2$ be closed oriented surfaces
with magnetic flows $\phi^{i}$ determined by pairs $(g_i,\Omega_i)$, $i=1,2$.
Suppose $\phi^1$ is $C^1$-conjugate to $\phi^2$ and let
$f:SM_1\to SM_2$ be the conjugacy. 
Then ${\frak m}(\mu^{1}_{\ell})={\frak m}(f_{*}\mu^{1}_{\ell})$.
In particular, $\phi^1$ has conjugate points if and only if $\phi^2$ does.

\label{l2}
\end{Lemma}

\begin{proof} Let $W_1(x,v)$ be the subspace of $T_{(x,v)}SM_1$ spanned by the magnetic vector
field $X_1$ and the vertical vector field $V$. 
Since $W_1$ contains the magnetic vector field and it is 2-dimensional it must be a 
Lagrangian subbundle. Since $df$ maps $X_1$ to $X_2$, it also maps Lagrangian subspaces
contained in $T(SM_1)$ to Lagrangian subspaces contained in $T(SM_2)$.
In particular, the subbundle $W_2:=df(W_1)$ must be a Lagrangian subbundle contained in $T(SM_2)$.

We now invoke the fact \cite{CGIP}[Corollary 3.2] that the asymptotic Maslov index of a
measure with zero asymptotic cycle does not depend on the continuous Lagrangian section 
that is used to compute it. Thus we can compute ${\frak m}(\mu^{1}_{\ell})$ using $W_1$
and ${\frak m}(f_{*}\mu^{1}_{\ell})$ using $W_2$ to readily obtain the 
equality claimed in the lemma.
To see that $\phi^1$ has conjugate points if and only if $\phi^2$ does we use
\cite{CGIP}[Theorem 4.4] which says that the asymptotic Maslov index of an invariant
probability measure (with zero asymptotic cycle) is positive if and only if there are 
conjugate points in its support.

\end{proof}

Let us now prove Theorem A. We know that the horocycle flow of a closed hyperbolic surface
has ${\frak a}(\mu_{\ell})={\frak m}(\mu_{\ell})=0$. Since the horocycle flow is transitive,
by Lemmas \ref{l1} and \ref{l2} the magnetic flow of $(g,\Omega)$ also has
${\frak a}(\mu_{\ell})={\frak m}(\mu_{\ell})=0$. 
The Proposition tells us that $g$ has constant curvature $k$ and $\Omega=\la\Omega_a$ 
with $k+\la^2=0$. If we let $a:=-k$, then $ag$ has curvature $-1$ and we
have the situation of two closed hyperbolic surfaces with $C^1$-conjugate
horocycle flows. The work of Marcus or Ratner \cite{M,Ra1}, tell us that $ag$ and $\bar{g}$ 
must in fact be isometric as desired.

\section{Closed Orbits in nontrivial free homotopy classes}

We consider magnetic flows defined on an arbitrary closed connected
manifold $M$. Let $g$ be a Riemannian metric and let $\Omega$ be a closed 2-form.
We will assume that $\Omega$ is {\it weakly exact}, that is,
the lift of $\Omega$ to $\M$, the universal covering of $M$, is exact.
Let $\theta$ be a primitive and let $c=c(g,\Omega)$ be Ma\~n\'e's critical value,
defined as in the Introduction. Recall that $c$ is finite if and only if
there is a bounded primitive.
If $\theta$ is a bounded primitive, then
our Lagrangian $L$ satisfies all the hypotheses of Aubry-Mather theory
for non compact manifolds as described for example in \cite{Co,Fa,FaMa}.
Recall that the {\it energy} in this case is simply the real valued function on $T\M$
given by $(x,v)\mapsto \frac{1}{2}|v|^{2}_{x}$.

We consider $\pi_{1}(M)$ acting on $\M$ by covering
tranformations and we let $\Pi:\M\to M$ be the covering projection.
Given a non-trivial element $\varphi\in\pi_{1}(M)$, let 
$Z_{\varphi}:=\{\rho\in\pi_{1}(M):\;\;\rho^{-1}\varphi\rho=\varphi\}$ be the
centralizer of $\varphi$.

\begin{Theorem} Let $k>c$ be given. Suppose that there exists a primitive
$\theta$ which is $Z_{\varphi}$-invariant and for which
\[\sup_{x\in\M}\frac{1}{2}\left|\theta_{x}\right|^2\leq k-\varepsilon\]
for some $\varepsilon>0$.
Then the non-trivial free homotopy class determined by $\varphi$ contains
a closed magnetic geodesic with energy $k$.
\label{previoatodo}
\end{Theorem}

\begin{proof}

Let $\theta$ be a $Z_{\varphi}$-invariant primitive with
\[\sup_{x\in\M}\frac{1}{2}\left|\theta_{x}\right|^2\leq k-\varepsilon\]
for some $\varepsilon>0$.
If we consider the Lagrangian on $\M$ given by
\[L(x,v)=\frac{1}{2}|v|^{2}_{x}-\theta_{x}(v)\]
then
\begin{equation}
L(x,v)+k\geq\varepsilon>0
\label{porabajo}
\end{equation}
for all $(x,v)\in TM$. Let us consider {\it Ma\~n\'e's action potential}, which is given by
\[\Phi_{k}(x,y)=\inf_{T>0}\Phi_{k}(x,y;T)\]
where
\[\Phi_{k}(x,y;T):=\inf_{\ga}A_{L+k}(\ga),\]
and $\ga$ ranges among all absolutely continuous curves defined on $[0,T]$
connecting $x$ to $y$. 
The potential $\Phi_k$ is a Lipschitz function which satisfies a triangle inequality
\[\Phi_{k}(x,y)\leq \Phi_{k}(x,z)+\Phi_{k}(z,y).\]
Note that the action potential is 
$Z_{\varphi}$-invariant, since $\theta$ is $Z_{\varphi}$-invariant.

Let $\psi\in\pi_{1}(M)$ be an arbitrary covering transformation. Since $\psi^*\theta-\theta$
is closed, there exists a smooth function $f_{\psi}:\M\to \re$ such that
$\psi^*\theta-\theta=df_{\psi}$. The function $f_{\psi}$ is uniquely defined up to addition
of a constant, so from now on we shall assume that $f_{\psi}$ is the unique function
for which $f_{\psi}(x_0)=0$ where $x_0$ is some fixed point in $\M$. Note that from the
definition of $\Phi_k$ we have:
\begin{equation}
\Phi_{k}(\psi x,\psi y)=\Phi_{k}(x,y)+f_{\psi}(y)-f_{\psi}(x)
\label{picara1}
\end{equation}
for all $x,y\in\M$ and all $\psi\in\pi_{1}(M)$.



\begin{Lemma} If $\psi_{1}^{-1}\varphi\psi_{1}=\psi_{2}^{-1}\varphi\psi_{2}$, 
then $f_{\psi_{1}}=f_{\psi_{2}}$.
\label{facil}
\end{Lemma}

\begin{proof} Clearly $\tau:=\psi_{1}\psi_{2}^{-1}\in Z_{\varphi}$.
Hence $\psi^{*}_{1}\theta-\theta=\psi_{2}^{*}\tau^{*}\theta-\theta=
\psi^{*}_{2}\theta-\theta$ which implies $df_{\psi_{1}}=df_{\psi_{2}}$.
Thus $f_{\psi_{1}}=f_{\psi_{2}}$ since they both vanish at $x_0$.

\end{proof}

A theorem due to Ma\~n\'e \cite{Ma,CDI} ensures that given two distinct points $x$ and $y$ in
$\M$ there exists a magnetic geodesic $\ga:[0,R]\to \M$ with energy $k$, which connects $x$ to $y$
and realizes the potential, i.e.,
\[A_{L+k}(\ga)=\Phi_{k}(x,y).\]
On account of (\ref{porabajo})
\begin{equation}
\Phi_{k}(x,y)\geq \varepsilon\,R\geq \frac{\varepsilon}{\sqrt{2k}}\,d(x,y).
\label{cotaabajo}
\end{equation}

Let $a:=\inf_{x\in\M}\Phi_{k}(x,\varphi x)$. Take a sequence of points $x_n$ such that
$\Phi_{k}(x_{n},\varphi x_{n})\to a$. Let $K$ be a compact fundamental domain for the action
of $\pi_{1}(M)$ on $\M$ and let $\psi_{n}\in\pi_{1}(M)$ be such that $\psi_{n}^{-1}x_{n}\in K$.
Let $y_{n}:=\psi^{-1}_{n}x_{n}$. Without loss of generality we can assume that $y_n$ converges
to some point $y\in K$.
Using the triangle inequality for $\Phi_{k}$ we have:
\[\Phi_{k}(\psi_{n}y,\varphi\psi_{n}y)\leq 
\Phi_{k}(\psi_{n}y,\psi_{n}y_{n})+\Phi_{k}(\psi_{n}y_{n},\varphi\psi_{n}y_{n})
+\Phi_{k}(\varphi\psi_{n}y_{n},\varphi\psi_{n}y).\] 
Using the $\varphi$-invariance of $\Phi_{k}$ and (\ref{picara1}) we obtain
\begin{align*}
\Phi_{k}(\psi_{n}y,&\varphi\psi_{n}y) \\
&\leq\Phi_{k}(y,y_{n})+f_{\psi_{n}}(y_{n})-f_{\psi_{n}}(y)
+\Phi_{k}(y_n,y)+f_{\psi_{n}}(y)-f_{\psi_{n}}(y_{n})+\Phi_{k}(x_{n},\varphi x_{n})\\
&=\Phi_{k}(y,y_{n})+\Phi_{k}(y_n,y)+\Phi_{k}(x_{n},\varphi x_{n}).
\end{align*}
But the expression $\Phi_{k}(y,y_{n})+\Phi_{k}(y_n,y)+\Phi_{k}(x_{n},\varphi x_{n})$ is bounded
in $n$, hence there exists $C>0$ such that
\[\Phi_{k}(\psi_{n}y,\varphi\psi_{n}y)<C\]
for all $n$. Inequality (\ref{cotaabajo}) now implies that there exist only finitely many different
elements of the form $\psi_{n}^{-1}\varphi\psi_{n}$. Hence for infinitely many values of $n$,
$\psi_{n}^{-1}\varphi\psi_{n}$ equals some fixed covering tranformation, let us say, $\lambda$. 
Without loss of generality we shall assume that $\psi_{n}^{-1}\varphi\psi_{n}=\lambda$ for all $n$.
Lemma \ref{facil} tells us that $f_{\psi_{n}}$ is independent of $n$, so let us
set $f_{0}:=f_{\psi_{n}}$.

Using (\ref{picara1}) again we have
\[\Phi_{k}(\psi_{n}y,\varphi\psi_{n}y)
=\Phi_{k}(y,\la y)+f_{0}(\la y)-f_{0}(y)\]
and
\[\Phi_{k}(x_{n},\varphi x_{n})=\Phi_{k}(\psi_{n}y_{n},\varphi\psi_{n}y_{n})
=\Phi_{k}(y_{n},\la y_{n})+f_{0}(\la y_n)-f_{0}(y_n).\]
Since $\Phi_{k}(x_n,\varphi x_n)\to a$ we conclude that
$\Phi_{k}(\psi_{n}y,\varphi\psi_{n}y)=a$ for all $n$, hence the
points $\psi_{n} y$ realize the infimum of the function
$x\mapsto \Phi_{k}(x,\varphi x)$.

Let $z$ be one of these points, i.e. $\Phi_{k}(z,\varphi z)=a$.
Consider the minimizing
magnetic geodesic $\ga$ with energy $k$ given by Ma\~n\'e's theorem which connects
$z$ to $\varphi z$ and for which
\[A_{L+k}(\ga)=\Phi_{k}(z,\varphi z).\]
We claim that $d\varphi(\dot{\ga}(0))=\dot{\ga}(R)$. This implies that the projection
of $\ga$ to $M$ gives a closed magnetic geodesic in the free homotopy class determined by
$\varphi$. To prove that $d\varphi(\dot{\ga}(0))=\dot{\ga}(R)$ we play the same game as in
Riemannian geometry. Consider $b>0$ small and note that
\begin{align*}
\Phi_{k}(\ga(b),\varphi\ga(b))&\leq \Phi_{k}(\ga(b),\varphi z)
+\Phi_{k}(\varphi z, \varphi\ga(b))\\
&= \Phi_{k}(\ga(b),\varphi z)+\Phi_{k}(z,\ga(b))\\
&= \Phi_{k}(z,\varphi z)
\end{align*}
where in the first equality we used the $\varphi$-invariance of $\Phi_{k}$.
Since $x\mapsto \Phi_{k}(x,\varphi x)$ achieves its minimum at $z$, 
we must have $d\varphi(\dot{\ga}(0))=\dot{\ga}(R)$ which concludes the proof of the theorem.

\end{proof}

The next lemma will be important for us. Its proof is a fairly standard application
of amenability.

\begin{Lemma} Let $\Gamma\subset\pi_{1}(M)$ be an amenable subgroup. For any $k>c$, there
exists a smooth $\Gamma$-invariant primitive $\vartheta$ such that
\[\sup_{x\in\M}\frac{1}{2}\left|\vartheta_{x}\right|^2\leq k.\]
\label{amenable}
\end{Lemma}

\begin{proof} Since $\Gamma$ is amenable it has a right invariant mean on $\ell^{\infty}(\Gamma)$, that is,
there exists a bounded linear functional $m:\ell^{\infty}(\Gamma)\to\re$ such that
\begin{enumerate}
\item $m(a)=a$ for a constant function $a$;
\item $m(a_1)\geq m(a_2)$, if $a_1(\varphi)\geq a_{2}(\varphi)$ for all $\varphi\in\Gamma$;
\item $m(\varphi_{*}a)=m(a)$, where $\varphi_{*}a(\psi):=a(\psi\varphi)$.
\end{enumerate}

By the definition of Ma\~n\'e's critical value, given $k>c$, there exists a primitive $\theta$
such that
\begin{equation}
\sup_{x\in\M}\frac{1}{2}\left|\theta_{x}\right|^2\leq k.
\label{bond}
\end{equation}
Given $(x,v)\in TM$, consider that function $a_{(x,v)}:\Gamma\to\re$ given by
$a_{(x,v)}(\varphi)=\theta_{\varphi(x)}(d\varphi(v))$.
Since $\varphi$ acts by isometries, inequality (\ref{bond}) implies
that $a_{(x,v)}\in \ell^{\infty}(\Gamma)$. Hence we can set:
\[\vartheta_{x}(v):=m(a_{(x,v)}).\]
The linearity and continuity of $m$ implies that $\vartheta$ is a smooth 1-form and by property (3), $\vartheta$
is $\Gamma$-invariant. Moreover if $\ga$ is any closed curve, the linearity and continuity of $m$ imply
\begin{equation}
\int_{\gamma}\vartheta=m\left(\varphi\mapsto \int_{\ga}\varphi^{*}\theta\right).
\label{rescue}
\end{equation}
But $\int_{\ga}\varphi^{*}\theta$ is independent of $\varphi$. This can be seen as follows.
Since $\M$ is simply connected, there exists a smooth map $F:{\mathbb D}\to\M$, where ${\mathbb D}$
is a 2-disk and the restriction of $F$ to $\partial{\mathbb D}$ is $\ga$.
Hence
\[\int_{\ga}\varphi^{*}\theta=\int_{\partial{\mathbb D}}F^{*}\varphi^*\theta
=\int_{{\mathbb D}}F^*\varphi^{*} d\theta.\]
But $d\theta=\widetilde{\Omega}$ and hence is $\Gamma$-invariant, so the last integral is independent
of $\varphi$. Thus from (\ref{rescue})
\[\int_{\ga}\vartheta=\int_{\ga}\theta\]
for any closed curve $\ga$. This implies that $\vartheta-\theta$ is an exact form.

Finally, since $a_{(x,v)}(\varphi)\leq \sqrt{2k}$ for all $x\in\M$, all $v\in T_{x}\M$ with norm one
and all $\varphi\in \Gamma$,
property (2) implies that $\vartheta_{x}(v)\leq \sqrt{2k}$ for all $x\in\M$ and all $v\in T_{x}\M$ with norm one.
Thus
\[\sup_{x\in\M}\frac{1}{2}\left|\vartheta_{x}\right|^2\leq k\]
as desired.

\end{proof}

\begin{Corollary}Suppose that $\pi_{1}(M)$ is amenable and $\Omega$ is not exact.
Then $c(g,\Omega)=\infty$.
\label{cinfinito}
\end{Corollary}

\begin{proof} If $c(g,\Omega)$ is finite, $\widetilde{\Omega}$ admits a bounded
primitive and by the previous lemma, $\widetilde{\Omega}$ admits a $\pi_{1}(M)$-invariant
primitive $\theta$. The form $\theta$ descends to $M$ showing that $\Omega$
is exact.
\end{proof}

\begin{Theorem} Let $k>c$ and let $\varphi\in\pi_{1}(M)$ be a non-trivial element
with amenable centralizer. Then the non-trivial free homotopy class determined by $\varphi$ 
contains a closed magnetic geodesic with energy $k$.
\label{cloorb}
\end{Theorem}

\begin{proof} It follows right away from Theorem \ref{previoatodo}
and Lemma \ref{amenable}.
\end{proof}

\begin{Remark}{\rm If $\pi_{1}(M)$ is the fundamental group of a closed
manifold of negative curvature and $\varphi$ is non-trivial, Preissman's
theorem implies that $Z_{\varphi}$ coincides with the infinite cyclic group 
generated by $\varphi$, which is of course amenable.
Thus we can apply Theorem \ref{cloorb} to any non-trivial free homotopy class.
In the next section we will apply the theorem to a closed surface of genus $\geq 2$.}
\end{Remark}

\subsection{Proof of Theorem C} We will need the following result
which was proven by G. Contreras \cite{Co1} in the exact case and extended
by O. Osuna \cite{O} to the weakly exact case as part of his Ph.D thesis work.

\begin{Theorem} For almost every $k$ in the interval $(0,c)$ there exists a
closed contractible magnetic geodesic with energy $k$.
\label{gonzalo-osvaldo}
\end{Theorem}

The first item in Theorem $C$ follows from Corollary \ref{cinfinito}
and Theorem \ref{gonzalo-osvaldo}.
The second item in Theorem C follows  from Theorem \ref{cloorb}
and Theorem \ref{gonzalo-osvaldo}.

Theorem \ref{gonzalo-osvaldo} is proved by showing that an appropriate action functional (which
depends on the energy level) on the space of contractible loops exhibits a mountain pass geometry. 
Then standard Morse theory gives the existence of critical points whenever the Palais-Smale condition
holds. It has been known for some time, that the Palais-Smale condition can only fail in the
``time direction". An argument originally due to M. Struwe can now be applied to the mountain pass
geometry to overcome this difficulty for almost every energy level.

\section{Proof of Theorem B}

In this section we return to the case in which $M$ is a closed surface of genus
$\geq 2$.
The following lemma has independent interest.

\begin{Lemma} Suppose the magnetic flow of $(g,\Omega)$ is Ma\~n\'e critical.
If there exists a nontrivial free homotopy class without closed magnetic geodesics, 
then there exists an invariant Borel probability measure
$\nu$ with ${\mathfrak a}(\nu)={\mathcal S}(\nu)=0$. Equivalently, the magnetic flow
is not of contact type. 
\label{contacto}
\end{Lemma}

\begin{proof} Let $\sigma$ be the nontrivial free homotopy class without closed magnetic 
geodesics. Suppose that $\sigma$ is generated by the covering transformation $\varphi$.
As in the proof of Theorem \ref{cloorb} we consider a $\varphi$-invariant action
potential $\Phi_{k}\geq 0$ for all $k> c=1/2$.  Now take a decreasing sequence $k_n$ approaching $1/2$ as $n\to\infty$.
Theorem \ref{cloorb} gives points $x_n$ and orbits $\ga_n:[0,T_n]\to\M$ with energy $k_n$
connecting $x_n$ and $\varphi x_n$. The orbits $\ga_n$ project to $M$ as closed orbits
in the class $\sigma$ and
\[0\leq A_{L_n+k_n}(\ga_n)=\Phi_{k_n}(x_n,\varphi x_n).\]
Moreover, $x_n$ is a minimum of $x\mapsto \Phi_{k_{n}}(x,\varphi x)$. Hence if $y$ is any
point in $\M$, 
\[\Phi_{k_n}(x_n,\varphi x_n)\leq \Phi_{k_{n}}(y,\varphi y)\leq C\]
for some constant $C>0$, since $L_n(x,v)+k_n\leq 1/2+k_{n}+\sqrt{2k_{n}}$ for all
$(x,v)\in TM$ with $|v|_{x}\leq 1$. Hence
\begin{equation}
0\leq A_{L_{n}+k_n}(\ga_n)\leq C.
\label{bound}
\end{equation}
We now observe that $\inf_{n} T_{n}>0$, otherwise we would get curves in the class
$\sigma$ with arbitrarily short lengths, which is impossible.
If $\sup_{n}T_{n}<\infty$, by passing to a subsequence if necessary, we can assume
that $T_{n}\to T_{0}$ and that the projection of $(\ga_{n}(0),\dot{\ga}_{n}(0))$
to $TM$ converges to some point $(p,v)\in SM$. The orbit of $(p,v)$ gives rise to
a closed magnetic geodesic with period $T_0$ in the homotopy class $\sigma$.
Since we are assuming that $\sigma$ has no such orbits we must have $\sup_{n}T_{n}=\infty$.
Without loss of generality we shall assume from now on that $T_n\to\infty$.

Let us indicate with a tilde the lift of objects on $M$ (or $SM$) to $\M$ (or $S\M$).
Note that
\begin{equation}
(L+k_n)(\ga_n,\dot{\ga}_{n})=2k_n-\theta_{n}(\dot{\ga}_{n})=
(\widetilde{\alpha}-\widetilde{\pi}^*\theta_n)(\widetilde{X})(\ga_{n},\dot{\ga}_{n}).
\label{eq1}
\end{equation}
Let $\Theta$ be a primitive of $\omega$ in a neighbourhood of $SM$.
Since 
$$d(\widetilde{\alpha}-\widetilde{\pi}^*\theta_n)=-\widetilde{\omega}=-d\widetilde{\Theta}$$
there exists a smooth closed 1-form $\rho_n$ defined in a neighbourhood of $S\M$ for which
\begin{equation}
\widetilde{\alpha}-\widetilde{\pi}^*\theta_n=-\widetilde{\Theta}+\rho_n.
\label{eq2}
\end{equation}
Combining (\ref{eq1}) and (\ref{eq2}) we obtain:
\begin{equation}
(L+k_n)(\ga_n,\dot{\ga}_{n})=
-\widetilde{\Theta}(\widetilde{X})(\ga_n,\dot{\ga}_{n})
+\rho_{n}(\widetilde{X})(\ga_n,\dot{\ga}_{n}).
\label{eq3}
\end{equation}
Let $\nu_n$ be the Borel probability measures on $TM$ given by:
\[\int f\,d\nu_n:=\frac{1}{T_{n}}\int_{0}^{T_{n}}f(\Pi(\ga_{n}(t)),d\Pi(\dot{\ga}_{n}(t)))\,dt.\]
Without loss of generality we can assume that $\nu_n$ converges weakly to
an invariant measure $\nu$. Since $k_n\to 1/2$, the measure $\nu$ has support
in $SM$. Let us check that ${\mathfrak a}(\nu)={\mathcal S}(\nu)=0$.
Using (\ref{bound}) and (\ref{eq3}) we see that
\begin{align*}
0=\lim \frac{1}{T_{n}}A_{L+k_n}(\ga_n)&=
-\int\Theta(X)\,d\nu+\lim \frac{1}{T_{n}}\int_{(\ga_n,\dot{\ga}_{n})}\rho_n\\
\end{align*}
and therefore to show that ${\mathfrak a}(\nu)=0$ it suffices to check that
\[\lim \frac{1}{T_{n}}\int_{(\ga_n,\dot{\ga}_{n})}\rho_n=0.\]
Equality (\ref{eq2}) implies that $\rho_n$ is $\varphi$-invariant and its norm
is bounded by a constant, let us say $A$, independent of $n$.
Let $\widehat{M}$ be the manifold obtained by taking the quotient of $\M$ by
the action of the cyclic group generated by $\varphi$. The curves $\ga_n$ project
to {\it simple} closed curves in $\widehat{M}$ which are all homotopic and therefore
the curves $(\ga_n,\dot{\ga}_{n})$ in $T\M$ project to closed curves $\Gamma_n$ in 
a neighbourhood of $S\widehat{M}$, whose homology class $[\Gamma_{n}]$ is independent
of $n$. The form $\rho_n$ descends to a closed 1-form $\widehat{\rho}_{n}$ defined 
in a neighbourhood of $S\widehat{M}$. Observe that
\[\int_{(\ga_n,\dot{\ga}_{n})}\rho_n=\int_{\Gamma_{n}}\widehat{\rho}_n=
\langle [\widehat{\rho}_{n}],[\Gamma_{n}]\rangle.\]
Since $[\Gamma_n]$ is independent of $n$ and $\widehat{\rho}_{n}$ is bounded by $A$
we have
\[\lim \frac{1}{T_n}\langle [\widehat{\rho}_{n}],[\Gamma_{n}]\rangle=0\]
as desired.

Let us prove that ${\mathcal S}(\nu)=0$.
Let $\Upsilon $ be any closed 1-form on $SM$. Since $\pi^*:H^{1}(SM,\re)\to H^{1}(M,\re)$
is an isomorphism, there exists a closed 1-form $\delta$ on $M$ and a smooth
function $G$ on $SM$ such that $\Upsilon=\pi^*\delta+dG$. Thus
\[\int\Upsilon(X)\,d\nu=\int\pi^*\delta(X)\,d\nu\]
and to prove that ${\mathcal S}(\nu)=0$ it suffices to show that
\[\int \pi^*\delta(X)\,d\nu=0.\]
But
\[\int \pi^*\delta(X)\,d\nu=\lim\frac{1}{T_{n}}\int_{\Pi\circ\ga_n}\delta=
\lim\frac{1}{T_{n}}\langle [\delta],[\Pi\circ\ga_{n}]\rangle .\]
The curves $\Pi\circ\ga_n$ are all in the same free homotopy class $\sigma$, hence
the homology class $[\Pi\circ\ga_n]$ is independent of $n$ which gives ${\mathcal S}(\nu)=0$
as desired.

To complete the proof of the lemma, recall that $SM$ is of contact type 
if and only if for all invariant Borel probability measures $\mu$ with zero
asymptotic cycle, ${\mathfrak a}(\mu)\neq 0$ (cf. \cite{CMP}[Proposition 2.4]).

\end{proof}

Let us prove Theorem B. The first observation is that we always have semistatic curves
starting at any point in $\M$ \cite{CIPP,Co,FaMa}.
These are magnetic geodesics $\ga:[0,\infty)\to \M$ such that
\[A_{L+1/2}(\ga|_{[s,t]})=\Phi_{1/2}(\ga(s),\ga(t))\] 
for $0\leq s<t<\infty$ where as before
\[\Phi_{1/2}(x,y)=\inf_{T>0}\Phi_{1/2}(x,y;T).\]
A semistatic curve must be free of conjugate points in $[0,\infty)$
(see \cite{CI}[Corollary 4.2]) and hence the $\omega$-limit set of the projection
of $\ga$ to $M$ must also be free of conjugate points.
Since we are assuming that the magnetic flow is uniquely ergodic, this implies
that all $SM$ is free of conjugate points.
On account of Lemma \ref{contacto} and unique ergodicity, ${\mathfrak a}(\mu_{\ell})=0$
and the theorem follows from the Proposition.

\begin{Remark}{\rm  We can rephrase Lemma \ref{contacto} by saying that if a Ma\~n\'e critical
magnetic flow is of contact type, then every nontrivial free homotopy class contains a closed
magnetic geodesic. Since there are nontrivial free homotopy classes with the property that any 
closed curve in them is homologous to zero, we obtain, in particular, closed magnetic geodesics
homologous to zero.
Recall that the {\it Weinstein conjecture} says that every Reeb vector field on a closed
3-manifold admits a closed orbit. The {\it strong Weinstein conjecture} asserts that in fact one can
find finitely many closed orbits which form a cycle homologous to zero.
There has been lots of progress recently regarding this conjecture (cf. \cite{ACH}).
However, the work of J. Entyre \cite{E} implies that magnetic flows, with $\Omega$ symplectic, 
are {\it excluded} from all the known cases in which the conjecture holds.
}

\label{enthof}
\end{Remark}

\end{document}